\documentclass[12pt,a4paper]{amsart}
\usepackage[utf8]{inputenc}	
\usepackage[T1]{fontenc}
\usepackage[in,headings]{fullpage}
\usepackage{amsbsy, amscd, amsfonts, amsmath, amsrefs, amssymb, amsthm}
\usepackage{graphicx}
\usepackage{float}
\usepackage{color}   
\usepackage[colorlinks=true,linkcolor=blue,citecolor=blue,urlcolor=blue]{hyperref}
\usepackage{comment}
\excludecomment{A}

\newtheorem{theorem}{Theorem}
\newtheorem*{theor}{Theorem}

\newtheorem{corollary}[theorem]{Corollary}
\newtheorem{proposition}[theorem]{Proposition}
\newtheorem{lemma}[theorem]{Lemma}

\newtheorem*{nonumberthm}{Theorem}

\theoremstyle{definition}

\newtheorem{remark}[theorem]{Remark}

\theoremstyle{remark}

\newtheorem{claim}{Claim}

\newcommand{\C}{\mathbf{C}}
\newcommand{\Z}{\mathbf{Z}}

\newcommand{\R}{\mathbf{R}}
\newcommand{\N}{\mathbf{N}}

\renewcommand{\Re}{\mathop{\mathrm{Re}}\nolimits}

\newfont{\cmbsy}{cmbsy10}
\newfont{\cmmib}{cmmib10}

\begin{document}

\title{Explicit van der Corput's $d$-th derivative estimate. }
\author[Arias de Reyna]{J. Arias de Reyna}
\address{%
Universidad de Sevilla \\ 
Facultad de Matem\'aticas \\ 
c/Tarfia, sn \\ 
41012-Sevilla \\ 
Spain.} 




\email{arias@us.es, ariasdereyna1947@gmail.com}


\begin{abstract}
We give an explicit version for van der Corput's $d$-th derivative estimate of exponential sums. 
\begin{theor} 
Let $X$, and $Y\in\R$ be such that $\lfloor Y\rfloor>d$ where $d\ge3$ is a natural number. 
Let $f\colon(X,X+Y]\to\R$ be a real function with continuous derivatives up to the order $d$. Assume that $0<\lambda\le f^{(d)}(x)\le\Lambda$ for $X<x\le X+Y$.  Denote by $D=2^d$. Then 
\begin{equation}
\Bigl|\frac{1}{Y}\sum_{X<n\le X+Y}e(f(n))\Bigr|\le
\max\Bigl\{A_d\Bigl(\frac{\Lambda}{\lambda Y}\Bigr)^{2/D}, B_d\Bigl(\frac{\Lambda^2}{\lambda}\Bigr)^{1/(D-2)},C_d(\lambda Y^d)^{-2/D}\Bigr\},
\end{equation}
where 
$A_d$, $B_d$, and $C_d$ are explicit constants. They depend on $d$ but for $d\ge2$ for example $A_d< 7.5$, $B_d<5.8$ and $C_d<10.9$.
\end{theor}

We follow the reasoning of van der Corput in three papers published in 1937, that contained an error. I correct this error and try to get the smallest possible constants.
We apply this theorem to zeta sums, giving the best choice of $d$ in each case. Also, we prove that our Theorem implies Titchmarsh's Theorem 5.13.

\end{abstract}

\maketitle
 
\section{Introduction}
The van der Corput $d$-th derivative estimate gives a bound of an exponential sum 
\[\Bigl|\frac{1}{Y}\sum_{n=X}^{X+Y}e^{2\pi i f(n)}\Bigr|,\]
assuming that $f\colon[X,X+Y]$ has continuous derivatives up to order $d$ and 
also satisfies the inequalities $0<\lambda\le f^{(d)}(x)\le \Lambda$. Here, we prove  an explicit version 
\begin{nonumberthm}
Let $X$, and $Y\in\R$ be such that $\lfloor Y\rfloor>d$, where $d\ge3$ is a natural number. 
Let $f\colon(X,X+Y]\to\R$ be a real function with continuous derivatives up to the order $d$. Assume that $0<\lambda\le f^{(d)}(x)\le\Lambda$ for $X<x\le X+Y$.  Denote by $D=2^d$. Then 
\begin{equation}\label{E:nonumber}
\Bigl|\frac{1}{Y}\sum_{X<n\le X+Y}e(f(n))\Bigr|\le
\max\Bigl\{A_d\Bigl(\frac{\Lambda}{\lambda Y}\Bigr)^{2/D}, B_d\Bigl(\frac{\Lambda^2}{\lambda}\Bigr)^{1/(D-2)},C_d(\lambda Y^d)^{-2/D}\Bigr\},
\end{equation}
where the constants $A_d$, $B_d$ and $C_d$ are given in Table \ref{table1}, all of them
less than 11.
\end{nonumberthm}
We do not find explicit versions in modern references. Titchmarsh's Theorem 5.13 is the best known. Our corollary \ref{TitchmarshTh} shows an explicit version of Titchmarsh's theorem as a consequence of our fundamental theorem. The Graham and Kolesnik book Theorem \cite{GK}*{Th.~2.8} is also non-explicit and has slightly worse exponents than our theorem.  Other non-explicit expositions of the van der Corput theorem can be found in Montgomery \cite{Mont}, Bordellès \cite{B}, Robert \cite{R}, and many others.

The techniques of van der Corput were refined by Vinogradov et al. Notably the recent results of Wooley \cite{W}, Bourgain et al. \cite{BDG}. These results have allowed considerable improvement of the exponents in theorems of van der Corput type, see, for example, the theorem proved by Heath-Brown \cite{HB}.

Our main interest in the subject is the search for explicit results. To this end, modern results are useless. It is true that in number theory, interest is generally focused on the order of magnitude, in which case modern theorems are much better. But in certain numerical studies, explicit results are needed. For example, to bound the error when truncating certain series or to improve the time cost of certain algorithms. We refer to the Hiary presentation \cite{H} for more details of this. 

The original result of van der Corput \cite{C}*{1928} was explicit. It contained instead of the supremum, the sum of the three terms on the right-hand side of \eqref{E:nonumber} all with constant 21. In \cite{C1}*{1937} and \cite{C3}*{1937} a version of the theorem was given, where it contained the maximum of the three terms, all with constant 25.  Shortly after the appearance of \cite{C} Titchmarsh \cite{T1}*{1931} published a simpler proof of a similar theorem (Theorem 5.13 in \cite{T}), not explicit but having practically the same utility. After this, the explicit versions of van der Corput were somewhat forgotten.

Since I needed to make use of an explicit version, my attention fell on van der Corput's second proof.  It was surprising that neither Titchmarsh \cite{T}, Montgomery \cite{Mont} nor Graham and Kolesnik \cite{GK} cited this work of van der Corput.
I found a problem in van der Corput's proof of the analog of our Theorem \ref{T:main}, which appears in \cite{C1}*{p.~670--671}.  The problem is in equality between line $-7$ and line $-6$ on page 671. This equality can be summarized in \[\frac{1}{N}\sum_{n=1}^{N}\max(A, a_n)=\max\Bigl(A,\frac{1}{N}\sum_{n=1}^N a_n\Bigr).\]
Since the function $f(x)=\max(A,x)$ is convex,  the correct sign would be $\ge$. The proof fails at this point. But the demonstration is easy to correct (see our Theorem \ref{T:main}).
The proof of Theorem \ref{T:main} is by induction, in our proof the main induction step is changed. I consider it useful to write it in detail.

After I wrote this some other explicit version of the $d$-th estimate have been published. 
See \cite{Y}.

\subsection{Some notation} We denote by $\N$ the set of natural numbers, excluding $0$. $A\subset B$ means $x\in A\Longrightarrow x\in B$, hence including the case $A=B$. For typographical reasons, we use the abbreviation $e(x)=e^{2\pi i x}$.  We will use the letter $S$ to denote a sum of exponentials and $E$ to denote the corresponding mean value, when $S$ is divided by the length of the range.

\section{Kusmin Landau lemma}

\begin{lemma}
Let $f\colon(X,X+Y]\to\R$ a real function continuously differentiable with a monotonous derivative  and such that $\theta\le f'(x)\le 1-\theta$ for some $0<\theta\le \frac12$, then
\[\Bigl|\sum_{X<n\le X+Y}e^{2\pi i f(n)}\Bigr|\le \cot\frac{\pi\theta}{2}.\]
\end{lemma}

\begin{proof}
This is proved in Landau \cite{L}. As this sharp version is not usually given and the original is in German, a detailed proof has been included in the TeX source of this document. To make it appear in the resulting pdf, you should interchange the \%\ symbols in lines 11 and 12 of the TeX file.
\end{proof}

\begin{remark}\label{R:intplace}
The sum does not change if we replace $f(x)$ by $f(x)+kx$ with $k\in\Z$, therefore the hypothesis on $f'$  can be substituted by  $\theta\le f'(x)-k\le 1-\theta$ for some fixed integer $k$.
\end{remark}
\begin{remark}
van der Corput noticed that there is a bound of the trigonometric sum depending only on $\theta$.
Kusmin gives a simple proof with a bound of type $A/\theta$. It was Landau \cite{L} who gave the sharp constant $\cot\frac{\pi\theta}{2}$, and proved that it is the best possible bound.
\end{remark}

\begin{A}
\begin{proof}[Proof of Lemma 1]
Assume first that $f'$ is not decreasing. 
Let $a$ and $b$ be the first and last integers in the range $(X, X+Y]$. If there is none or $b=a$, the lemma is trivially true, for $b=a+1$ 
we have
\begin{align*}
|e^{2\pi i f(a)}+e^{2\pi i f(a+1)}|&=|(e^{2\pi i f(a)}+e^{2\pi i f(a+1)})e^{-\pi i(f(a)+f(a+1))}|\\
&=|e^{-\pi i (f(a+1)-f(a))}+e^{\pi i (f(a+1)-f(a))}|\\
&=2|\cos\pi(f(a+1)-f(a))|\\
&\le \cos\pi\theta+\cos\pi\theta\le \frac{\cos\pi\theta}{\sin\pi\theta}+\frac{1}{\sin\pi\theta}=\cot\frac{\pi\theta}{2}.
\end{align*}
In other cases, we prepare the sum to apply partial summation.
\begin{align*}
\sum_{X<n\le X+Y}&e^{2\pi i f(n)}=\\
&=\sum_{n=a}^{b-1}\frac{e^{2\pi i f(n)}-e^{2\pi i f(n+1)}}{1-e^{2\pi i(f(n+1)-f(n))}}+e^{2\pi i f(b)}=\sum_{n=a}^{b-1}(e^{2\pi i f(n)}-e^{2\pi i f(n+1)})F(n)+e^{2\pi i f(b)}\\
&=e^{2\pi i f(a)}F(a)+\sum_{n={a+1}}^{b-1} e^{2\pi i f(n)}(F(n)-F(n-1))-e^{2\pi i f(b)}F(b-1)+e^{2\pi i f(b)}
\end{align*}
So, 
\[\Bigl|\sum_{X<n\le X+Y}e^{2\pi i f(n)}\Bigr|\le |F(a)|+\sum_{n=a+1}^{b-1}|F(n)-F(n-1)|+|1-F(b-1)|\]
Let $g(x)=f(x+1)-f(x)$, then 
\[F(n)=\frac{1}{1-e^{2\pi i(f(n+1)-f(n))}}=\frac12+\frac12\frac{e^{-\pi i g(n)}+e^{\pi i g(n)}}{e^{-\pi i g(n)}-e^{\pi i g(n)}}=\frac12+\frac{i}{2}\cot(\pi g(n)).\]
Hence,
\[|F(n)-F(n-1)|=\frac12|\cot(\pi g(n))-\cot(\pi g(n-1))|=\frac12\cot(\pi g(n-1))-\frac12\cot(\pi g(n)),\]
The last equality is true because $g(n)=f(n+1)-f(n)=f'(\xi_1)\ge f'(\xi_2)=g(n-1)$ since $f'$ is not decreasing, and $\cot(x)$ is a decreasing function.
The sum collapses, and we obtain
\[
\Bigl|\sum_{X<n\le X+Y}e^{2\pi i f(n)}\Bigr|\le |F(a)|+\frac12\cot(\pi g(a))-\frac12\cot(\pi g(b-1))+|1-F(b-1)|\]
Notice that 
\[|F(a)|^2=\frac{1+\cot^2(\pi g(a))}{4}=\frac{1}{4\sin^2(\pi g(a))},\quad
\text{and}\quad |1-F(b-1)|^2=\frac{1}{4\sin^2(\pi g(b-1))}.\]
Therefore,
\begin{align*}
\Bigl|\sum_{X<n\le X+Y}e^{2\pi i f(n)}\Bigr| &\le \frac{1}{2\sin(\pi g(a))}+\frac12\cot(\pi g(a))-\frac12\cot(\pi g(b-1))+
\frac{1}{2\sin(\pi g(b-1))}\\
&= \frac{1+\cos(\pi g(a))}{2\sin(\pi g(a))}+\frac{1-\cos(\pi g(b-1))}{2\sin(\pi g(b-1))}\le 
\frac{1+\cos(\pi \theta)}{\sin(\pi \theta)}=\cot\frac{\pi\theta}{2}.
\end{align*}
The changes in the case $f'$ is nonincreasing affect only the sign of 
the expression for $|F(n)-F(n-1)|$, without consequences. 
\end{proof}

\end{A}

\begin{lemma}\label{L:2ChGr}
Let $f\colon(X,X+Y]\to\R$ a real valued function two times differentiable with continuity. Assume that $0<\lambda\le f''(x)\le \Lambda$, then 
\[\Bigl|\sum_{X<n\le X+Y}e(f(n))\Bigr|\le (\Lambda Y+2)\Bigl(1+\frac{4}{\sqrt{\pi\lambda}}\Bigr)+1.\]
\end{lemma}
\begin{proof}
Since $f'$ is continuous, the image $f'(X,X+Y]$ is an interval $I$. For any two points in $I$, we have $|f'(b)-f'(a)|\le |f''(\xi)(b-a)|\le \Lambda Y$. Let $[a,b]$ be the least closed interval with extremes $a$, $b\in\Z$ and such that $[a,b]\supset I$. The bound on the length of $I$ implies 
$0\le b-a\le \Lambda Y+2$. 
For any $0<\theta<1/2$ we split the sum 
\begin{align*}
\sum_{X<n\le X+Y}e(f(n))&=\sum_{a\le f'(n)\le a+\theta}e(f(n))+
\sum_{k=a+1}^{b-1}\sum_{k-\theta<f'(n)\le k+\theta}e(f(n))+\sum_{b-\theta <f'(n)\le b}e(f(n))\\&+
\sum_{k=a}^{b-1}\sum_{k+\theta<f'(n)\le k+1-\theta}e(f(n)).
\end{align*}
If there are two integers $n$ and $m$ with $a\le f'(n),f'(m)\le a+\theta$, we have \[\theta \ge |f'(n)-f'(m)|= |(n-m)f''(\xi)|\ge |n-m|\lambda.\] Therefore, the first sum contains $0$, $1$, or $\le (1+\theta/\lambda)$ terms. In any case $\le (1+\theta/\lambda)$ terms. The same happens for the sum on $b-\theta <f'(n)\le b$. 
A similar reasoning gives that each sum on $k-\theta<f'(n)\le k+\theta$ contains at most 
$(1+2\theta/\lambda)$.  By the trivial bound, all these sums are bounded by 
\[(1+\theta/\lambda)+(b-a-1)(1+2\theta/\lambda)+(1+\theta/\lambda)=
(b-a)(1+2\theta/\lambda)+1\le (\Lambda Y+2)(1+2\theta/\lambda)+1.\] 

To each sum $\sum_{k+\theta<f'(n)\le k+1-\theta}e(f(n))=\sum_{k+\theta<f'(n)\le k+1-\theta}e(f(n)-kn)$ we may apply Kusmin-Landau's lemma (see remark \ref{R:intplace}). 
Note that $f'(x)-k$ is monotonous. 
Therefore, 
\begin{align*}
\Bigl|\sum_{X<n\le X+Y}e(f(n))\Bigr|&\le (\Lambda Y+2)(1+2\theta/\lambda)+1+(b-a)\cot\frac{\pi\theta}{2}\\
&\le (\Lambda Y+2)\Bigl(1+2\theta/\lambda+\cot\frac{\pi\theta}{2}\Bigr)+1.
\end{align*}
Since $\cot\alpha\le\frac{1}{\alpha}$
\begin{equation}\label{E:partial}
\Bigl|\sum_{X<n\le X+Y}e(f(n))\Bigr|\le (\Lambda Y+2)\Bigl(1+2\theta/\lambda+\frac{2}{\pi\theta}\Bigr)+1.
\end{equation}
Taking $\theta=\sqrt{\lambda/\pi}$ we obtain 
\[\Bigl|\sum_{X<n\le X+Y}e(f(n))\Bigr|\le (\Lambda Y+2)\Bigl(1+\frac{4}{\sqrt{\pi\lambda}}\Bigr)+1.\]
This will end the proof if $\theta\le 1/2$ that is, when $\lambda\le \pi/4$, When $\lambda>\pi/4$, the sum is bounded by the trivial bound as $Y+1$, therefore, 
\[\frac{1}{\Lambda Y+2}\Bigl|\sum_{X<n\le X+Y}e(f(n))\Bigr|\le\frac{Y+1}{\Lambda Y+2}.\]
as a function of $Y\in[0,\infty)$, $\frac{Y+1}{\Lambda Y+2}$ is monotonous and varies between 
$\frac12$ and $\frac{1}{\Lambda}$, both $\le \bigl(1+\frac{4}{\sqrt{\pi\lambda}}\bigr)$ in fact 
$\frac12<1$ and $\frac{1}{\Lambda}\le\frac{1}{\lambda}$. This is less than $\frac{4}{\sqrt{\pi\lambda}}$ for $\lambda>\pi/16$.
\end{proof}

The next lemma, due to van der Corput, is a version of Lemma \ref{L:2ChGr}. 
\begin{lemma}\label{C:maincor}
Let $f\colon(X,X+Y]\to\R$ a real-valued function twice differentiable with continuity. Assume that $0<\lambda\le f''(x)\le \Lambda$, and $Y\ge 1$, then 
\[\Bigl|\sum_{X<n\le X+Y}e(f(n))\Bigr|\le \frac{A}{\sqrt{\lambda}}(\Lambda Y+2),\]
with $A=\frac{2}{\sqrt{\pi}}(1+\sqrt{1+3 \pi/8})=2.79368380731\dots$
\end{lemma}
\begin{proof}
We assume the hypothesis of Lemma \ref{L:2ChGr}, therefore we have \eqref{E:partial}.
Assume also that $\lambda\le\beta^2$ for some $\beta>0$ to be determined, then 
the sum $S:=\sum e(f(n))$ satisfies,   taking $\theta=\alpha\sqrt{\lambda}$. If $\alpha\sqrt{\lambda}<1/2$, we have
\[\frac{|S|}{\Lambda Y+2}\le 1+\frac{2\theta}{\lambda}+\frac{2}{\pi\theta}+\frac{1}{\Lambda Y+2}\le 
\frac{1.5\beta}{\sqrt{\lambda}}+\frac{2\alpha}{\sqrt{\lambda}}+\frac{2}{\pi\alpha\sqrt{\lambda}}=\Bigl(1.5\beta+2\alpha+\frac{2}{\pi\alpha}\Bigr)\frac{1}{\sqrt{\lambda}}.\]
The best option here is to take $\alpha=1/\sqrt{\pi}$ so we look for a  constant $A$ such that 
\[1.5\beta+\frac{4}{\sqrt{\pi}}\le A.\]
For $\lambda<\beta^2$, and our final choice: $A=\frac{1}{\beta}=2.793\dots$, we have $\theta=\alpha\sqrt{\lambda}=\sqrt{\lambda/\pi}<\beta/\sqrt{\pi}<1/2$.

For $\lambda\ge\beta^2$ we use the trivial bound for the sum
\[\frac{|S|}{\Lambda Y+2}\le\frac{Y+1}{\Lambda Y+2}.\]
As a function of $Y$ for $Y\ge1$ this is a monotonous function taking values between 
$\frac{1}{\Lambda}\le \frac{1}{\lambda}$ and $\frac{2}{2+\Lambda}\le \frac{2}{2+\lambda}$. We want these two quantities to be bounded by $A/\sqrt{\lambda}$. In particular $1/A\le \sqrt{\lambda}$ must be true for 
$\lambda\ge\beta^2$. Therefore $1/A\le\beta$.  We also need that 
$\frac{2}{2+\lambda}\le A/\sqrt{\lambda}$ for $\lambda\ge\beta^2$.

We take $A$ as  $1/\beta$ where $\beta$ is defined by the equation 
\[\frac{1}{\beta}=1.5\beta+\frac{4}{\sqrt{\pi}}.\]
The solution gives us 
\[A=\frac{2}{\sqrt{\pi}}(1+\sqrt{1+3 \pi/8})=2.79368380731\dots\]
It is mechanical  to check that with this election we have 
$\frac{2}{2+\lambda}\le A/\sqrt{\lambda}$ for $\lambda\ge\beta^2$.
\end{proof}

\section{Weyl-van der Corput Lemma}

The idea of Weyl was to bound the square of the sum of exponentials to get a new easier sum of exponentials. van der Corput added a clever use of the Cauchy-Schwarz inequality to 
get a sum with fewer terms.  

\begin{lemma}[Weyl-van der Corput]\label{L:WvdC}
Let $f\colon(X, X+Y]\to\C$ a complex function with $|f(x)|\le 1$ and $H\le Y$ a natural number, then 
\begin{equation}\label{E:WvdC}
\Bigl|\frac{1}{Y}\sum_{X<n\le X+Y} f(n) \Bigr|^2\le \frac{4}{H}+\frac{4}{HY}\Re\Bigl(\frac{1}{H}\sum_{a=1}^{H-1}(H-a)\sum_{X<n,n+a\le X+Y} f(n+a)\overline{f(n)}\Bigr).
\end{equation}
\end{lemma}
\begin{proof}
It is convenient to extend $f$ putting $f(x)=0$ for $x\not\in (X, X+Y]$.  These values do not appear in equation \eqref{E:WvdC}. 

Notice that we have
\[H\sum_{X<n\le X+Y} f(n)=\sum_{X<n\le X+Y+H-1}\sum_{a=0}^{H-1}f(n-a)\]
Now,  Cauchy-Schwarz's inequality yields
\[H^2\Bigl|\sum_{X<n\le X+Y} f(n)\Bigr|^2\le (Y+H)\sum_{X<n\le X+Y+H-1}\Bigl|\sum_{a=0}^{H-1}f(n-a)\Bigr|^2.\]
Each individual square can be expanded
\[\Bigl|\sum_{a=0}^{H-1}f(n-a)\Bigr|^2=\sum_{a=0}^{H-1}f(n-a) \sum_{b=0}^{H-1}\overline{f(n-b)}.\]
We divide each of these sums into three parts
\[\sum_{a=b}+ \sum_{a<b}+\sum_{a>b} f(n-a)\overline{f(n-b)}, \qquad \text{with  }\quad  0\le a, b\le H-1.\]
In each of these sums a fixed  $n$ appears, some terms can be $=0$ due to the extension of $f$ we have made. In the sum with $a=b$, since $|f(x)|\le 1$,  we will have 
\[\sum_{a=b}f(n-a)\overline{f(n-b)}=\sum_{a=0}^{H-1}|f(n-a)|^2\le H.\]
The other two sums are complex conjugates, so their sum is twice their real part. Therefore,
\[H^2\Bigl|\sum_{X<n\le X+Y} f(n)\Bigr|^2\le H(Y+H)^2+2(Y+H)\mskip-15mu\sum_{X<n\le X+Y+H-1}\mskip-5mu \Re\Bigl(\sum_{0\le a<b\le H-1}
f(n-a)\overline{f(n-b)}\Bigr).\]
or  
\begin{align*}
H^2\Bigl|&\sum_{X<n\le X+Y} f(n)\Bigr|^2\\
&\le H(Y+H)^2+2(Y+H)\Re\Bigl(\sum_{X<n\le X+Y+H-1}\mskip10mu \sum_{0\le a<b\le H-1}
f(n-a)\overline{f(n-b)}\Bigr).
\end{align*}
For each term, we have $n-a>n-b$, so the non-null terms are of the form $f(m+k)\overline{f(m)}$, where $X<m<m+k\le Y$. 

Once fixed $m$ and $m+k$ with $X<m<m+k\le X+Y$, we can take $(a,b)$ with $0\le a<b\le H-1$ giving this product only for $b=a+k$. The possible pairs are as follows
\[(a,b)= (0,k), (1, k+1), (2,k+2),\dots, (H-k-1, H-1).\]
These are $(H-k)$ pairs. For each of these pairs, the corresponding value of $n=m+b$ satisfies the condition $X<n\le X+Y+H-1$.  It follows that the sum can be written as follows
\[\sum_{X<n\le X+Y+H-1}\sum_{0\le a<b\le H-1}
f(n-a)\overline{f(n-b)}\Bigr)=\sum_{m,k}(H-k) f(m+k)\overline{f(m)},\]
where the sum is in al $k$ with $1\le k\le H-1$ and only the terms with 
$X<m<m+k\le X+Y$ are not null.
Therefore,
\[H^2\Bigl|\sum_{X<n\le X+Y} f(n)\Bigr|^2\le H(Y+H)^2+2(Y+H)\Re\Bigl(\sum_{k=1}^{H-1}(H-k)\mskip-10mu \sum_{X<m<m+k\le X+Y} f(m+k)\overline{f(m)}\Bigr).\]
Since $H\le Y$ 
\[H^2\Bigl|\sum_{X<n\le X+Y} f(n)\Bigr|^2\le 4HY^2+4Y\Re\Bigl(\sum_{k=1}^{H-1}(H-k) \sum_{X<m<m+k\le X+Y}f(m+k)\overline{f(m)}\Bigr).\]
We get \eqref{E:WvdC} dividing by $H^2Y^2$. 
\end{proof}

\section{Induction and main result of van der Corput}

We need to introduce some notation to simplify the statement of the main result. We will start with a function $f\colon (X, X+Y]\to\R$, with $X$,  $Y\in\R$ and $Y>0$, and we get bounds for the sum of exponentials
\begin{equation}
E(f):=\frac{1}{Y}\sum_{X<n\le X+Y} e^{2\pi i f(x)}=\frac{1}{Y}\sum_{X<n\le X+Y} e(f(x)).
\end{equation}
For typographical reasons, we will use the standard notation $e(x)=e^{2\pi i x}$. 

By induction, we define new functions starting with $f$. Given  $\mathbf{a}=(a_1,a_2,\dots, a_h)$ with $a_j\in \R$, we define
\begin{align*}
f_{(a_1)}(x)&=f(x+a_1)-f(x),\\
f_{(a_1,a_2)}(x)&=f_{(a_1)}(x+a_2)-f_{(a_1)}(x),\\
\cdots & \cdots\\
f_{(a_1,a_2,\dots, a_h)}(x)&=f_{(a_1,a_2,\dots, a_{h-1})}(x+a_h)-f_{(a_1,a_2,\dots, a_{h-1})}(x).
\end{align*}
With $\mathbf{a}=(a_1,a_2,\dots, a_h)$, we put $\mathbf{a'}=(a_1,\dots, a_{h-1})$, and the last equation is equivalent to 
\[f_{\mathbf{a}}(x)=f_{\mathbf{a'}}(x+a_h)-f_{\mathbf{a'}}(x).\]
If $f$ is  $d$-times differentiable with continuity, we find by induction the integral representation
\begin{equation}\label{E:incInt}
f_{\mathbf{a}}(x)=\int_0^{a_d}\cdots\int_0^{a_2}\int_0^{a_1}f^{(d)}(x+t_d+\cdots+t_1)\,dt_1\,dt_2\,\cdots\,dt_d,\qquad \mathbf{a}=(a_1,a_2,\dots, a_d).
\end{equation}
This shows that $f_{\mathbf{a}}(x)$ does not depend on the order of elements in $\mathbf{a}$, but  this property does not depend on the differentiability of $f$.

In the next theorem we will use the constant 
\begin{equation}\label{E:DefB}
B=2+2\sqrt{2}=4.8284271247\dots,\qquad\text{with the property}\quad \frac{4}{H}+4T\le \max(B^2H^{-1}, BT).
\end{equation}
Because if $\frac{4}{H}\le (B-4)T$, then $\frac{4}{H}+4T\le BT$ in the other case, if $\frac{4}{H}> (B-4)T$ we have 
\[\frac{4}{H}+4T\le\frac{4}{H}+\frac{16}{(B-4)H}=\frac{4B}{B-4}\frac{1}{H}=\frac{B^2}{H}.\]

The domain of definition $f_{\mathbf{a}}$ is only a subset of $(X, X+Y]$ that may be empty.  But notice that $f_{\mathbf{a}}(x)$ is a linear  combination of $f(x+a_J)$ where $J\subset\{1,2,\dots, h\}$ and $a_J:=\sum_{j\in J} a_j$. We will use only positive $a_j$ and in this case $f_{\mathbf{a}}$ is defined exactly in $(X, X+Y-\sum_{j=1}^h a_j]$.
We will use the notation $I({\mathbf{a}})$ to denote this domain of the natural definition of $f_{\mathbf{a}}$.
However, we \emph{use the convention that $e(f_{\mathbf{a}}(x))=0$ if the value of $f_{\mathbf{a}}(x)$ is not defined.} This simplifies the range of our sums.
\begin{theorem}\label{T:main}
Let  $f\colon(X,X+Y]\to\R$ be a real function,  $d$, and $H_1$, \dots, $H_d$ natural numbers  with $H_1+H_2+\cdots+H_d\le Y$, then 
\begin{equation}\label{E:induction}
|E(f)|:=\Bigl|\frac{1}{Y}\sum_{X<n\le X+Y}e(f(x))\Bigr|\le B\max(H_1^{-1/2}, H_2^{-1/4},\dots, H_d^{-1/2^d}, B^{-1/2^d}T_d^{1/2^d}),
\end{equation}
where $B$ is the constant defined in \eqref{E:DefB} and 
\[T_d=\frac{1}{YH_1\cdots H_d}\sum_{\mathbf{a}}\Bigl|\sum_{n\in I(\mathbf{a})}e(f_{\mathbf{a}}(n))\Bigr|,\]
where $\textbf{a}$ run through all vectors $(a_1,a_2,\dots a_d)$ with integers  $a_r$ such that $1\le a_r\le H_r-1$.
\end{theorem}
\begin{proof}
We proceed by induction on $d$. For $d=1$, the Weyl-van der Corput  lemma \eqref{L:WvdC} yields 
\[|E(f)|^2\le \frac{4}{H_1}+\frac{4}{H_1Y}\Re\Bigl(\frac{1}{H_1}\sum_{a=1}^{H_1-1}(H_1-a)\sum_{X<n,n+a\le X+Y} e(f(n+a)-f(n))\Bigr).\]
This easily implies that 
\[|E(f)|^2\le \frac{4}{H_1}+\frac{4}{H_1Y}\sum_{a=1}^{H_1-1}\Bigl|\sum_{n\in I{(a)}} e(f_{(a)}(n))\Bigr|=\frac{4}{H_1}+4T_1\le \max\Bigl(\frac{B^2}{H_1}, B T_1\Bigr).\]
So,
\[|E(f)|\le \max\bigl(B H_1^{-1/2}, \sqrt{B}\,T_1^{1/2}\bigr)\le B\max\bigl(H_1^{-1/2}, B^{-1/2}T_1^{1/2}\bigr),\]
and \eqref{E:induction} is proved for $d=1$.

For $d>1$, by the induction hypothesis, since $H_1+\cdots H_{d-1}\le Y$ we have
\begin{equation}\label{E:191127-1}
|E(f)|\le B\max(H_1^{-1/2}, \dots,H_{d-1}^{-1/2^{d-1}}, B^{-1/2^{d-1}}T_{d-1}^{1/2^{d-1}}),
\end{equation}
where 
\[T_{d-1}=\frac{1}{YH_1\cdots H_{d-1}}\sum_{\mathbf{a'}}\Bigl|\sum_{n\in I(\mathbf{a'})}e(f_{\mathbf{a'}})(n))\Bigr|,\qquad \mathbf{a'}=(a_1,\dots, a_{d-1}).\]
The domain $I(\mathbf{a'})=(X, X+Y-(a_1+\cdots a_{d-1})]$. Therefore, we have 
$H_d\le Y-(H_1+\cdots+H_{d-1})\le Y-(a_1+\cdots a_{d-1})$ and we can apply Lemma \eqref{L:WvdC} to $e((f_{\mathbf{a'}})(n))$.
This yields
\begin{multline*}
\Bigl|\frac{1}{Y-(a_1+\cdots+a_{d-1})}\sum_{n\in I(\mathbf{a'})}e(f_{\mathbf{a'}})(n))\Bigr|^2\le \\
\frac{4}{H_d}+\frac{4}{H_d(Y-(a_1+\cdots+a_{d-1}))}\sum_{a_d=1}^{H_d-1}\Bigl|\sum_{n, n+a_d\in I(\mathbf{a'})} e(f_{\mathbf{a'}}(n+a_d)-f_{\mathbf{a'}}(n))\Bigr|.
\end{multline*}
By definition, $f_{\mathbf{a'}}(n+a_d)-f_{\mathbf{a'}}(n)=f_{\mathbf a}(n)$. Therefore, we have 
\begin{multline*}
\Bigl|\frac{1}{Y-(a_1+\cdots+a_{d-1})}\sum_{n\in I(\mathbf{a'})}e(f_{\mathbf{a'}})(n))\Bigr|^2
\\ \le 
\frac{4}{H_d}+\frac{4}{H_d(Y-(a_1+\cdots+a_{d-1}))}\sum_{a_d=1}^{H_d-1}\Bigl|\sum_{n\in I(\mathbf{a})}e(f_{\mathbf a}(n))\Bigr|.\end{multline*}
Multiplying by $(Y-(a_1+\cdots+a_{d-1}))^2$ we can change the factor on the right side   by $Y$, and dividing by $Y^2$ we get 
\[\Bigl|\frac{1}{Y}\sum_{n\in I(\mathbf{a'})}e(f_{\mathbf{a'}})(n))\Bigr|^2\le 
\frac{4}{H_d}+\frac{4}{H_dY}\sum_{a_d=1}^{H_d-1}\Bigl|\sum_{n\in I(\mathbf{a})}e(f_{\mathbf a}(n))\Bigr|.\]
Hence, by the Schwarz inequality
\begin{align*}
T_{d-1}^2&=\Bigl\{\frac{1}{YH_1\cdots H_{d-1}}\sum_{\mathbf{a'}}\Bigl|\sum_{n\in I(\mathbf{a'})}e(f_{\mathbf{a'}})(n))\Bigr|\Bigr\}^2\le \frac{1}{H_1\cdots H_{d-1}}\sum_{\mathbf{a'}}\Bigl|\frac1Y\sum_{n\in I(\mathbf{a'})}e(f_{\mathbf{a'}})(n))\Bigr|^2\\
&\le \frac{1}{H_1\cdots H_{d-1}}\sum_{\mathbf{a'}}\Bigl(\frac{4}{H_d}+\frac{4}{H_dY}\sum_{a_d=1}^{H_d-1}\Bigl|\sum_{n\in I(\mathbf{a})}e(f_{\mathbf a}(n))\Bigr|\Bigr)\\
&\le \frac{4}{H_d}+
\frac{4}{Y}\frac{1}{H_1\cdots H_d}\sum_{\mathbf{a'}}\sum_{a_d=1}^{H_d-1}\Bigl|\sum_{n\in I(\mathbf{a})}e(f_{\mathbf{a}})(n))\Bigr|.
\end{align*}
This can be written 
\[T_{d-1}^2\le \frac{4}{H_d}+
\frac{4}{Y}\frac{1}{H_1\cdots H_d}\sum_{\mathbf{a}}\Bigl|\sum_{n\in I(\mathbf{a})}e(f_{\mathbf{a}})(n))\Bigr|\]

Hence,
\[T_{d-1}^2\le\frac{4}{H_d}+
\frac{4}{Y}\frac{1}{H_1\cdots H_d}\sum_{\mathbf{a}}\Bigl|\sum_{n\in I(\mathbf{a})}e(f_{\mathbf{a}})(n))\Bigr|=
\frac{4}{H_d}+4T_d\le \max(B^2H_d^{-1}, BT_d).\]
Joining this with \eqref{E:191127-1} yields
\begin{align*}
|E(f)|&\le B\max(H_1^{-1/2}, \dots,H_{d-1}^{-1/2^{d-1}}, B^{-1/2^{d-1}}(\max(B^2H_d^{-1}, BT_d))^{1/2^d})\\
&\le B\max(H_1^{-1/2}, \dots,H_{d-1}^{-1/2^{d-1}}, H_{d}^{-1/2^{d}}, B^{-1/2^d}T_d^{1/2^d}).\qedhere
\end{align*}
\end{proof}

\section{Specific bounds for the sum}
To apply Theorem \ref{T:main} we need bounds of $T_d$ and also select adequate values of the $H_j$'s. In this section, we show how to achieve a practical solution to these problems.
The integral representation \eqref{E:incInt} is useful to 
bound $T_d$ in terms  of the derivatives of $f$. 

We will need two lemmas to minimize our bounds, choosing some variables adequately.

\begin{lemma}\label{L:quince}
Let $\delta$ be a natural number and $\xi$ and $Y$ real numbers such that $1\le\xi\le Y^\delta$. Then there exist positive real numbers $x_1$, $x_2$, \dots, $x_\delta$ such that
\begin{equation}
\prod_{n=1}^\delta x_n=\xi,\qquad \min(\xi^{2^{n}/(2^{\delta+1}-2)}, Y^{2^n/2^\delta})\le
x_n\le Y.
\end{equation}
\end{lemma}

\begin{proof}
Separate the proof in two cases depending on whether $\xi^{2^\delta}\le Y^{2^{\delta+1}-2}$ or not.

(a) If $\xi^{2^\delta}\le Y^{2^{\delta+1}-2}$, we take for $1\le n\le \delta$
\[x_n:=\xi^{2^n/(2^{\delta+1}-2)}\le \xi^{2^\delta/(2^{\delta+1}-2)}\le Y.\]
In this case, the product $\prod_{n=1}^\delta x_n=\xi$, because $\sum_{n=1}^\delta2^n=2^{\delta+1}-2$. 

(b) $\xi^{2^\delta}> Y^{2^{\delta+1}-2}$. Find $N$ the least natural number such that 
\begin{equation}\label{E:b1}
\xi^{2^{N+1}}>Y^{(\delta-N+1)2^{N+1}-2}.
\end{equation}
For $N=\delta-1$ this inequality is $\xi^{2^\delta}>Y^{2^{\delta+1}-2}$ that is true by hypothesis. 
For $N=0$, the inequality is $\xi^2>Y^{2\delta}$. This is false, since we assume $\xi\le Y^\delta$.
Therefore $1\le N\le \delta-1$ and the inequality must be false for $N-1$ so that 
\begin{equation}\label{E:b2}
\xi^{2^N}\le Y^{(\delta-N+2)2^N-2}.
\end{equation}
Then we take $x_n$ as 
\[x_n=\begin{cases} (\xi Y^{N-\delta})^{2^n/(2^{N+1}-2)} & 1\le n\le N,\\
Y & N+1\le n\le \delta. \end{cases}\] 
The  product is easily computed 
\[x_1\cdots x_\delta=(x_1\cdots x_N)(x_{N+1}\cdots x_\delta)=(\xi Y^{N-\delta})Y^{\delta-N}=\xi.\]
We now prove $x_n\le Y$. For $N+1\le n\le\delta$ this is trivial. For $1\le n\le N$ we have
\[x_n=(\xi Y^{N-\delta})^{2^n/(2^{N+1}-2)}\le (\xi Y^{N-\delta})^{2^N/(2^{N+1}-2)}\]
and the inequality $(\xi Y^{N-\delta})^{2^N/(2^{N+1}-2)}\le Y$ is equivalent to $\xi^{2^N}\le Y^{(\delta-N)2^N+2^{N+1}-2}$ and this is equivalent to \eqref{E:b2}.

We end the proof showing that $x_n\ge Y^{2^n/2^\delta}$. For $N+1\le n\le \delta$ we have
$x_n=Y\ge  Y^{2^n/2^\delta}$.  

In the case $1\le n\le N$, note that \eqref{E:b1} is equivalent to $\xi Y^{N-\delta}>Y^{1-2^{-N}}$ so that 
\[x_n=(\xi Y^{N-\delta})^{2^n/(2^{N+1}-2)}>Y^{2^n/2^{N+1}} \]
Since $\delta\ge N+1$ we have $2^n/2^{N+1}\ge 2^n/2^\delta$ and this ends the proof.
\end{proof}

\begin{lemma}\label{L:191103-3}
Let  $M$, $N$, $Z$, $\alpha$ and $\beta$ positive real numbers. There exists  $0<\xi\le Z$ such that
\begin{equation}
\max(\xi^{-\alpha}, M \xi^{-\beta}, N\xi^\beta)\le \max(Z^{-\alpha}, MZ^{-\beta}, M^{\frac12}N^{\frac12},N^{\frac{\alpha}{\alpha+\beta}}).
\end{equation}
\end{lemma}
\begin{proof}
Let us consider the function $h(\xi):=\max(\xi^{-\alpha}, M \xi^{-\beta}, N\xi^\beta)$.
There is a point $0<\xi_0\le Z$ such that $h(\xi_0)\le h(\xi)$ for $0<\xi\le Z$. This is true since  $\lim_{\xi\to0^+}h(\xi)=+\infty$. 

If $\xi_0=Z$, then  $N\xi_0^\beta\le \max(\xi_0^{-\alpha}, M \xi_0^{-\beta})$. In the other case $N\xi_0^\beta>\max(\xi_0^{-\alpha}, M \xi_0^{-\beta})$ implies  $h(\xi)\le h(\xi_0)$ para $\xi_0-\varepsilon<\xi<\xi_0$, contradicting the definition of  $\xi_0$. Therefore, in this case  
\[h(\xi_0)=\max(\xi_0^{-\alpha}, M \xi_0^{-\beta})=\max(Z^{-\alpha}, M Z^{-\beta})\le \max(Z^{-\alpha}, MZ^{-\beta}, M^{\frac12}N^{\frac12},N^{\frac{\alpha}{\alpha+\beta}}).\]

If $0<\xi_0<Z$, we will have $\xi_0^{-\alpha}=N\xi_0^\beta$, or, $M\xi_0^{-\beta}=N\xi_0^\beta$. In the other case, if the maximum were taken in only one of the terms, $h(\xi_0-\delta)$ or  $h(\xi_0+\delta)$ will be less than $h(\xi_0)$.  A similar situation occurs if the maximum is taken at a point $\xi_0$ where 
 $\xi_0^{-\alpha}=M\xi_0^{-\beta}>N\xi_0^\beta$.  

When  $\xi_0^{-\alpha}=N\xi_0^\beta$, $N^{-1}=\xi_0^{\alpha+\beta}$ and
\[h(\xi_0)=\xi_0^{-\alpha}=N^{\frac{\alpha}{\alpha+\beta}}\le \max(Z^{-\alpha}, MZ^{-\beta}, M^{\frac12}N^{\frac12},N^{\frac{\alpha}{\alpha+\beta}}).\]

When  $M\xi_0^{-\beta}=N\xi_0^\beta$, $\xi_0^{2\beta}=M/N$ and
\[h(\xi_0)=N\xi_0^\beta=N(M/N)^\frac12=(MN)^{\frac12}\le \max(Z^{-\alpha}, MZ^{-\beta}, M^{\frac12}N^{\frac12},N^{\frac{\alpha}{\alpha+\beta}}).\qedhere\]
\end{proof}

\begin{lemma}\label{L:cased2}
Let $X$, and $Y\in\R$ be such that $Y\ge1$. Let $f\colon(X,X+Y]\to\R$ be a real function with continuous derivatives up to the order $2$. Assume that $0<\lambda\le f''(x)\le\Lambda$  for $X<x\le X+Y$. Then, 
\[\Bigl|\frac{1}{Y}\sum_{X<n\le X+Y}e(f(n))\Bigr|\le A\Bigl\{\Bigl(\frac{\Lambda^2}{\lambda }\Bigr)^{1/2}+2(\lambda Y^2)^{-1/2}\Bigr\},\]
where $A$ is the constant defined in Lemma \ref{C:maincor}.
\end{lemma}
\begin{proof}
Lemma \ref{C:maincor} yields
\[\Bigl|\sum_{X<n\le X+Y}e(f(n))\Bigr|\le \frac{A}{\sqrt{\lambda}}(2+\Lambda Y).\]
Then we have 
\[
\frac{A}{\sqrt{\lambda}}(2+\Lambda Y)\le AY \Bigl(\frac{2}{Y\sqrt{\lambda}}+
\frac{\Lambda }{\sqrt{\lambda}}\Bigr)=AY\Bigl\{\Bigl(\frac{\Lambda^2}{\lambda }\Bigr)^{1/2}+2(\lambda Y^2)^{-1/2}\Bigr\}.\]
This ends the proof, but later we will need a slightly weaker result.
For any $r>0$ we have $U+V\le \max\{(1+r)U, (1+r^{-1})V\}$. Taking $r$ so that 
$(1+r)A=2^{1/4} B$, we get $r=2^{1/4}B/A-1=1.05535\dots$, and  $C_2:=2(1+r^{-1})A=10.881690878\dots$
Therefore, we have
\begin{equation}\label{E:case2}
\Bigl|\frac{1}{Y}\sum_{X<n\le X+Y}e(f(n))\Bigr|\le \max\Bigl\{B_2\Bigl(\frac{\Lambda^2}{\lambda }\Bigr)^{1/2},C_2(\lambda Y^2)^{-1/2}\Bigr\},
\end{equation}
with $B_2=2^{1/4}B=5.74199\dots$, and $C_2=\frac{2^{5/4}AB}{2^{1/4}B-A}=10.881690\dots$. 
\end{proof}

\begin{theorem} \label{T:dthtest}
Let $X$, and $Y\in\R$ be such that $\lfloor Y\rfloor>d$ where $d\ge3$ is a natural number. 
Let $f\colon(X,X+Y]\to\R$ be a real function with continuous derivatives up to the order $d$. Assume that $0<\lambda\le f^{(d)}(x)\le\Lambda$ for $X<x\le X+Y$.  Denote by $D=2^d$. Then 
\begin{equation}\label{E:mainmain}
\Bigl|\frac{1}{Y}\sum_{X<n\le X+Y}e(f(n))\Bigr|\le
\max\Bigl\{A_d\Bigl(\frac{\Lambda}{\lambda Y}\Bigr)^{2/D}, B_d\Bigl(\frac{\Lambda^2}{\lambda}\Bigr)^{1/(D-2)},C_d(\lambda Y^d)^{-2/D}\Bigr\},
\end{equation}
where 
\begin{equation}\begin{split}
A_d=B\sqrt{2}\;&(d-2)^{2/D},\quad B_d=B((3^4/2^3)(2/3)^{2d}2^{D/4}A^2B^{-2})^{1/(D-2)},\\
&C_d=B(AD/B)^{4/D}(d-2)^{2(d-2)/D}.
\end{split}\end{equation}
\end{theorem}
\begin{proof}

Denote by $\delta=d-2$. We apply Theorem \ref{T:main}. Let $H_1$, $H_2$, \dots, $H_\delta$ natural numbers with $H_1+\cdots+H_\delta\le Y$, then 
\begin{equation}\label{E:startineq}
\Bigl|\frac{1}{Y}\sum_{X<n\le X+Y}e(f(x))\Bigr|\le B\max(H_1^{-1/2}, H_2^{-1/4},\dots, H_\delta^{-1/2^\delta}, B^{-1/2^\delta}T_\delta^{1/2^\delta}),
\end{equation}
where
\[T_\delta=\frac{1}{H_1\cdots H_\delta}\sum_{a_1,\dots, a_\delta}\Bigl|\frac{1}{Y}\sum_{n\in I(\mathbf{a})}e(f_{\mathbf{a}}(n))\Bigr|\]
In this sum $1\le a_r\le H_r-1$  and $f_{\mathbf{a}}(n)$ is defined for $X<n\le X+Y-a_1-\cdots -a_\delta$.  Notice that 
\[Y\ge Y-a_1-\cdots -a_\delta\ge Y-(H_1-1)-\cdots-(H_\delta-1)\ge \delta\ge 1.\]
This allows one to apply Lemma \ref{C:maincor} to the function $f_{\mathbf{a}}$.

We have 
\[f_{\mathbf{a}}(x)=\int_0^{a_\delta}\cdots\int_0^{a_2}\int_0^{a_1}f^{(\delta)}(x+t_\delta+\cdots+t_1)\,dt_1\cdots\,dt_\delta.\]
Therefore (recall $\delta=d-2$)
\[f''_{\mathbf{a}}(x)=\int_0^{a_\delta}\cdots\int_0^{a_2}\int_0^{a_1}f^{(d)}(x+t_\delta+\cdots+t_1)\,dt_1\cdots\,dt_\delta>\lambda a_1a_2\cdots a_\delta,\]
and also $f''_{\mathbf{a}}(x)\le\Lambda a_1a_2\cdots a_\delta$. 
Hence, Lemma \ref{C:maincor} yields
\[\Bigl|\sum_{n\in I(\mathbf{a})}e(f_{\mathbf{a}}(n))\Bigr|\le \frac{A(2+ Y\Lambda a_1a_2\cdots a_\delta)}{\sqrt{\lambda a_1a_2\cdots a_\delta} }=\frac{A}{\sqrt{\lambda} }\Bigl(\frac{2}{\sqrt{a_1a_2\cdots a_\delta}}+Y\Lambda\sqrt{a_1a_2\cdots a_\delta}\Bigr).\]
Therefore, for any $H_1+H_2+\cdots+H_{\delta}\le Y$ we have
\[T_\delta\le \frac{1}{H_1\cdots H_\delta Y}\sum_{a_1,\dots, a_\delta}\frac{A}{\sqrt{\lambda} }\Bigl(\frac{2}{\sqrt{a_1a_2\cdots a_\delta}}+Y\Lambda\sqrt{a_1a_2\cdots a_\delta}\Bigr),\]
where each $a_r$ runs through $1\le a_r\le H_r-1$. Since 
\[\sum_{a=1}^{H-1}\frac{1}{\sqrt{a}}\le \int_0^{H-1}\frac{dt}{\sqrt{t}}=2\sqrt{H-1}\le \frac{2H}{\sqrt{H+1}},\quad\text{and}\quad
\sum_{a=1}^{H-1}\sqrt{a}\le \int_1^H\sqrt{t}\,dt\le \frac23 H^{3/2}.\]
It follows that 
\[T_\delta\le \frac{A}{Y \sqrt{\lambda}}\Bigl(\frac{2^{\delta+1}}{\sqrt{(H_1+1)\cdots(H_\delta+1)}}+(2/3)^\delta Y\Lambda\sqrt{H_1\cdots H_\delta}\Bigr)\]
\[T_\delta\le \frac{2A}{Y \sqrt{\lambda}}\max\Bigl(\frac{2^{\delta+1}}{\sqrt{(H_1+1)\cdots(H_\delta+1)}},(2/3)^\delta Y\Lambda\sqrt{H_1\cdots H_\delta}\Bigr).\]
In our main inequality \eqref{E:startineq} $T_\delta$ appear with exponent $1/2^\delta$ as $B^{-1/2^\delta}T^{1/2^\delta}$. Since 
$\delta=d-2$ and $D=2^d$, we have $2^\delta=D/4$.  Hence we are interested in 
\begin{align*}
\Bigl(\frac{2A}{BY \sqrt{\lambda}}\frac{2^{\delta+1}}{\sqrt{(H_1+1)\cdots(H_\delta+1)}}\Bigr)^{4/D}&=
\Bigl(\frac{4A^2}{B^2Y^2 \lambda}\frac{D^2/4}{(H_1+1)\cdots(H_\delta+1)}\Bigr)^{2/D}\\
&=\sqrt{2}\,M((H_1+1)\cdots(H_\delta+1))^{-2/D},
\end{align*}
with
\[\sqrt{2}\,M=(A^2 B^{-2} D^2 Y^{-2}\lambda^{-1})^{2/D}.\]
Analogously, we have 
\begin{multline*}
\Bigl(\frac{2A}{BY \sqrt{\lambda}}(2/3)^\delta Y\Lambda\sqrt{H_1\cdots H_\delta}\Bigr)^{4/D}=\Bigl(
\frac{4A^2}{B^2Y^2 \lambda}(2/3)^{2d-4}Y^2\Lambda^2 H_1\cdots H_\delta\Bigr)^{2/D}\\=\sqrt{2}\,N(H_1\cdots H_\delta)^{2/D}
\end{multline*}
where \[\sqrt{2}\,N=(3^4 2^{-2} A^2 B^{-2} (2/3)^{2d} \Lambda^2 \lambda^{-1})^{2/D}.\]
Hence,
\[B^{-1/2^\delta}T_\delta^{1/2^\delta}\le \max\bigl(\sqrt{2}\,M((H_1+1)\cdots(H_\delta+1))^{-2/D},
\sqrt{2}\,N(H_1\cdots H_\delta)^{2/D}\bigr),\]
and by \eqref{E:startineq} we get, for any natural numbers $H_1$, \dots, $H_\delta$
with $H_1+\cdots+H_\delta\le Y$
\begin{multline*}
\Bigl|\frac{1}{Y}\sum_{X<n\le X+Y}e(f(x))\Bigr|\\
\le B\max\bigl(H_1^{-1/2}, H_2^{-1/4},\dots, H_\delta^{-1/2^\delta}, \sqrt{2}\,M((H_1+1)\cdots(H_\delta+1))^{-\frac2D},
\sqrt{2}\,N(H_1\cdots H_\delta)^{\frac2D}\bigr)\\
\le B\sqrt{2}\,\max\bigl((H_1+1)^{-1/2}, \dots, (H_\delta+1)^{-1/2^\delta}, M((H_1+1)\cdots(H_\delta+1))^{-\frac2D},
N(H_1\cdots H_\delta)^{\frac2D}\bigr).\\
\end{multline*}
The objective of the last step will be seen in the proof of claim \ref{Claim1} below. 
In the next steps, we apply Lemma \ref{L:quince} and \ref{L:191103-3} to pick up adequate values of the $H_j$. 
We will change this inequality in several forms given in the following claims:

\begin{claim}\label{Claim1}
For any positive real numbers $x_1$, \dots, $x_\delta$ with $x_1+\cdots+x_\delta\le Y$ we have
\[\Bigl|\frac{1}{Y}\sum_{X<n\le X+Y}e(f(n))\Bigr|\le B\sqrt{2}\,\max(x_1^{-1/2},\dots, x_\delta^{-1/2^\delta},
M(x_1\cdots x_\delta)^{-2/D}, N(x_1\cdots x_\delta)^{2/D}).\]
\end{claim}
The first member of this inequality is $\le(Y+1)/Y\le  2\le B\sqrt{2}$, so that the inequality is trivially true in the case of some of the $x_k\le 1$. If all $x_k\ge 1$, take $H_k=\lfloor x_k\rfloor$, and Claim \ref{Claim1} follows because $x_k^{-1/2^k}> (H_k+1)^{-1/2^k}$, $x_k^{2/D}\ge H_k^{2/D}$  and the inequality shown above. 

\begin{claim}\label{Claim2}
For any $0\le \xi\le (Y/\delta)^\delta$ the inequality 
\[\Bigl|\frac{1}{Y}\sum_{X<n\le X+Y}e(f(n))\Bigr|\le B\sqrt{2}\,\max((Y/\delta)^{-4/D}, \xi^{-2/(D-4)}, M\xi^{-2/D}, N\xi^{2/D})\]
holds.
\end{claim}
The inequality is true for $0<\xi\le1$ since in this case $B\sqrt{2}\,\xi^{-2/(D-4)}\ge2$ and the left-hand side is bounded by 2. Therefore, we may assume $1\le \xi\le (Y/\delta)^\delta$. Then by Lemma \ref{L:quince} there exist positive real numbers $x_n$ for $1\le n\le \delta$ such that 
\[\prod_{n=1}^\delta x_n=\xi;\qquad \min(\xi^{{2^n}/{(2^{\delta+1}-2)}}, (Y/\delta)^{{2^n}/{2^{\delta}}})\le x_n\le Y/\delta, \qquad 1\le n\le \delta.\]
Therefore, these numbers $x_n$ satisfies  $x_1+\cdots+x_\delta\le Y$ and  
\[x_n^{-1/2^n}\le \max(\xi^{-2/(D-4)}, (Y/\delta)^{-4/D}).\]
Combining these inequalities with the  Claim \ref{Claim1} instance for these numbers $x_n$ we get Claim \ref{Claim2}.

\begin{claim}\label{Claim3}
The inequality 
\[\Bigl|\frac{1}{Y}\sum_{X<n\le X+Y}e(f(n))\Bigr|\le B\sqrt{2}\,\max((Y/\delta)^{-4/D},  
M(Y/\delta)^{-2\delta/D}, M^{1/2}N^{1/2},  N^{D/(2D-4)})\]
holds.
\end{claim}
The lemma \ref{L:191103-3} shows that there is some $\xi\in(0,(Y/\delta)^\delta]$ such that 
\[\max(\xi^{-2/(D-4)}, M\xi^{-2/D}, N\xi^{2/D})\le 
\max((Y/\delta)^{-2\delta/(D-4)}, M(Y/\delta)^{-2\delta/D}, M^{\frac12}N^{\frac12},N^{\frac{D}{2D-4}}).\]
Combining this with claim \ref{Claim2} yields
\begin{multline*}
\Bigl|\frac{1}{Y}\sum_{X<n\le X+Y}e(f(n))\Bigr|\\
\le B\sqrt{2}\max((Y/\delta)^{-4/D},  (Y/\delta)^{-2\delta/(D-4)}, M(Y/\delta)^{-2\delta/D}, M^{1/2}N^{1/2},  N^{D/(2D-4)}).
\end{multline*}
Since $\delta=d-2$ and $D=2^d$, we have $4/D\le 2\delta/(D-4)$ for $d\ge3$, and since $Y> d>\delta$
we have 
\[(Y/\delta)^{-2\delta/(D-4)}\le  (Y/\delta)^{-4/D}.\]
This proves claim \ref{Claim3}.

Now we bound each of the terms in the inequality in Claim \ref{Claim3}.
\begin{align*}
B\sqrt{2}(MN&)^{1/2}=B(\sqrt{2}M\sqrt{2}N)^{1/2}\\&=B(A^2 B^{-2} D^2 Y^{-2}\lambda^{-1})^{1/D}
(3^4 2^{-2} A^2 B^{-2} (2/3)^{2d} \Lambda^2 \lambda^{-1})^{1/D}\\
&=B(3^2 2^{-1}(2/3)^d A^2B^{-2} D)^{2/D}\Bigl(\frac{\Lambda }{\lambda Y}\Bigr)^{2/D}= A'_d\Bigl(\frac{\Lambda}{\lambda Y}\Bigr)^{2/D},\quad\text{for $d\ge3$.}
\end{align*}
Since $\lambda\le \Lambda$ and $Y\ge  d>\delta$,  we have
\begin{align*}
B\sqrt{2}(Y/\delta)^{-4/D}&\le B\sqrt{2}(Y/\delta)^{-2/D}\le  B\sqrt{2}(Y\lambda/\delta\Lambda)^{-2/D}\\&= B\sqrt{2}\delta^{2/D} \Bigl(\frac{\Lambda}{\lambda Y}\Bigr)^{2/D}=A_d\Bigl(\frac{\Lambda}{\lambda Y}\Bigr)^{2/D}.
\end{align*}
\begin{align*}
B\sqrt{2}\cdot N^{D/(2D-4)}&\le B\sqrt{2}\cdot 2^{-D/(4D-8)}(3^4 2^{-2} A^2 B^{-2} (2/3)^{2d} \Lambda^2 \lambda^{-1})^{1/(D-2)}\\
&=B\cdot 2^{\frac{2d-3}{D-2}+\frac{D}{4(D-2)}}\cdot 3^{\frac{4-2d}{D-2}}\cdot (A/B)^{\frac{2}{D-2}}\Bigl(\frac{\Lambda^2}{\lambda}\Bigr)^{1/(D-2)}\\
&= B((81/8)(2/3)^{2d}2^{D/4}A^2B^{-2})^{1/(D-2)}\Bigl(\frac{\Lambda^2}{\lambda}\Bigr)^{1/(D-2)}=B_d\Bigl(\frac{\Lambda^2}{\lambda}\Bigr)^{1/(D-2)}.
\end{align*}
\begin{align*}
B\sqrt{2} M(Y/\delta)^{-2\delta/D}&=B(A^2 B^{-2} D^2 Y^{-2}\lambda^{-1})^{2/D}(Y/\delta)^{-2\delta/D}\\&
=B(AD/B)^{4/D}\delta^{2\delta/D}\cdot(Y^{2+\delta}\lambda)^{-2/D}\\
&=B(AD/B)^{4/D}\delta^{2\delta/D}\cdot(Y^{d}\lambda)^{-2/D}=C_d (\lambda Y^d)^{-2/D},\quad\text{for $d\ge3$.}
\end{align*}
It is easy to show that $A'_d\le A_d$ and equation \eqref{E:mainmain} is proved. 
\end{proof}

\begin{remark}
The sequences of the coefficients $A_d$, $B_d$, and $C_d$ are convergent. 
\[\lim_{d\to\infty}A_d=B\sqrt{2},\quad \lim_{d\to\infty}B_d=B 2^{1/4},\quad 
\lim_{d\to\infty}C_d=B.\]
\end{remark}

\begin{remark}
The numerical upper bounds for the coefficients  are given in the next table
\begin{table}[htp]
\caption{Coefficients in \eqref{E:mainmain}}
\begin{center}
\begin{tabular}{cllr}
$d$ & $A_d$  & $B_d$   & $C_d\ \ \ $ \\ \hline
2 &     0     &  5.742  & 10.882 \\
3 & 6.829 & 4.971 & 10.389 \\
4 & 7.447 & 5.094 & 10.016 \\
5 & 7.314 & 5.286 & 8.545 \\
6 & 7.131 & 5.445 & 7.197\\
7 & 7.003 & 5.558 & 6.264\\
8 & 6.925 & 5.632 & 5.679\\
9 & 6.881 & 5.678 & 5.324\\
10 & 6.857 & 5.706 & 5.114\\
11 & 6.844 & 5.722 & 4.992\\
$>21$  & 6.829 & 5.742 & 4.829
\end{tabular}
\end{center}
\label{table1}
\end{table}%
\end{remark}
\begin{remark}
Taking the supremum of the coefficients, we obtain an expression valid for all $d\ge 3$
\begin{equation}
\Bigl|\frac{1}{Y}\sum_{X<n\le X+Y}e(f(n))\Bigr|\le 
B\max\Bigl\{2^{5/8}\Bigl(\frac{\Lambda}{\lambda Y}\Bigr)^{2/D}, 2^{1/4}\Bigl(\frac{\Lambda^2}{\lambda}\Bigr)^{1/(D-2)},2\sqrt{2A/B}(\lambda Y^d)^{-2/D}\Bigr\}.
\end{equation}
And for all $d\ge2$, under the hypothesis of Theorem \ref{T:dthtest} or Lemma \ref{L:cased2}
\begin{equation}
\Bigl|\frac{1}{Y}\sum_{X<n\le X+Y}e(f(n))\Bigr|\le 
\max\Bigl\{7.447\Bigl(\frac{\Lambda}{\lambda Y}\Bigr)^{2/D}, 5.742\Bigl(\frac{\Lambda^2}{\lambda}\Bigr)^{1/(D-2)},10.882(\lambda Y^d)^{-2/D}\Bigr\}.
\end{equation}
\end{remark}

\section{Applications}

In the Theory of the Zeta function we are interested in the zeta-sums
\[\sum_{n\le X}\frac{1}{n^\sigma}e^{-it\log n}=\overline{\sum_{n\le X} \frac{1}{n^\sigma} e(f(n))}, \qquad 
f(x)=\frac{t}{2\pi}\log x.\]
By partial summation, we may reduce to the case $\sigma=0$.  Usually the range of $n$ is separated in intervals, so that in the end the objective is to bound 
\begin{equation}
S(X,t):=\sup_{0<Y\le X}\Bigl|\sum_{X<n\le X+Y} e(f(n))\Bigr|,\qquad X=\Bigl(\frac{t}{2\pi}\Bigr)^\alpha,
\end{equation}
where $\alpha$ will be called the exponent of the sum $S(X,t)$, they are convenient to express the bounds. 

\begin{proposition}\label{P:30}
Let $X>0$ and $t>2\pi$ be real numbers and $d\ge 2$ be an integer such that $X^{1-2/D}\ge d$ with $D=2^d$, then 
\begin{equation}\label{E:coef32}
S(X,t)\le \tau^{\alpha}\max\Bigl\{\widehat{A_d} \,\tau^{-2\alpha/D}, \widehat{B_d}\, \tau^{\frac{1-\alpha d}{D-2}},\widehat{C_d}\, \tau^{-2/D}\Bigr\}, \qquad \tau=t/2\pi, \quad X=\tau^\alpha,
\end{equation}
where the constants $\widehat{A_d}$, $\widehat{B_d}$ and $\widehat{C_d}$ are given in \eqref{constants}
\end{proposition}
\begin{proof}
Let $d<Y\le X$, we apply Theorem \ref{T:dthtest} (or Lemma \ref{L:cased2} \eqref{E:case2} for $d=2$). The derivatives of our function 
$f(x)=\frac{t}{2\pi}\log x$ are 
\[f'(x)=\frac{t}{2\pi x}, \quad f''(x)=-\frac{t}{2\pi x^2},\qquad f^{(k)}(x)=(-1)^{k+1}\frac{t\; (k-1)!}{2\pi x^k}.\]
Therefore, for $X<x\le X+Y\le 2X$ we have 
\[\lambda = \frac{t(d-1)!}{2\pi 2^d X^d} \le  |f^{(d)}(x)|\le \frac{t(d-1)!}{2\pi X^d} =\Lambda\]
It follows that 
\[\Bigl|\sum_{X<n\le X+Y}e(f(n))\Bigr|\le Y\max\Bigl\{\widehat{A_d} Y^{-2/D},\widehat{B_d} \Bigl(\frac{t}{2\pi X^d}\Bigr)^{1/(D-2)},\widehat{C_d}\Bigl(\frac{t Y^d}{2\pi X^d}\Bigr)^{-2/D}\Bigr\},\]
where
\begin{equation}\label{constants}
\begin{aligned}
\widehat{A_d}&=2^{2d/D}A_d\\
\widehat{B_d}&=((d-1)!\, D)^{1/(D-2)} B_d \\
\widehat{C_d}&=((d-1)!/D)^{-2/D} C_d
\end{aligned}
\end{equation}
In the main inequality, the exponent of $Y$ in the three terms are $\ge0$ so that we can change $Y$ by its limit superior $X$. 
\[\Bigl|\sum_{X<n\le X+Y}e(f(n))\Bigr|\le X\max\Bigl\{\widehat{A_d} X^{-2/D},\widehat{B_d} \Bigl(\frac{t}{2\pi X^d}\Bigr)^{1/(D-2)},\widehat{C_d}\Bigl(\frac{t }{2\pi }\Bigr)^{-2/D}\Bigr\}.\]
When $d\ge3$, since $ \widehat{A_d} X^{1-2/D}\ge 6d$, this inequality is also true for $Y\le d$, for $d=2$ the coefficient $\widehat{A_2}=0$, but the second term is $\ge \widehat{B_2}>11$ and the inequality is also true for $Y\le 2$. Therefore, taking sup in $Y$ we obtain
\[S(X,t)\le X\max\Bigl\{\widehat{A_d} X^{-2/D},\widehat{B_d} \Bigl(\frac{t}{2\pi X^d}\Bigr)^{1/(D-2)},\widehat{C_d}\Bigl(\frac{t }{2\pi }\Bigr)^{-2/D}\Bigr\}.\]
Now by hypothesis $X=\tau^{\alpha}$ with $\tau=t/2\pi$ so that
\[S(X,t)\le \tau^{\alpha}\max\Bigl\{\widehat{A_d} \,\tau^{-2\alpha/D}, \widehat{B_d}\, \tau^{\frac{1-\alpha d}{D-2}},\widehat{C_d}\, \tau^{-2/D}\Bigr\}.\qedhere\]
\end{proof}
\begin{remark}
Upper bounds of the new coefficients $\widehat{A_d}$, $\widehat{B_d}$ and $\widehat{C_d}$, are given in  table \ref{table2}.
\begin{table}[htp]
\caption{Coefficients in \eqref{E:coef32}}
\begin{center}
\begin{tabular}{crrr}
$d$ & $\widehat{A_d}$  & $\widehat{B_d}$   & $\widehat{C_d}\ \ \ $ \\ \hline
2 &     0     &  11.484  & 21.764 \\
3 & 11.484 & 7.891 & 14.692 \\
4 & 10.531 & 7.058 & 11.323 \\
5 & 9.083 & 6.596 & 8.700 \\
6 & 8.121 & 6.290 & 7.057\\
7 & 7.554 & 6.086 & 6.098\\
8 & 7.232 & 5.953 & 5.548\\
9 & 7.051 & 5.869 & 5.234\\
10 & 6.950 & 5.817 & 5.056\\
11 & 6.895 & 5.786 & 4.955\\
$>21$  & 6.829 & 5.743 & 4.829
\end{tabular}
\end{center}
\label{table2}
\end{table}%
\end{remark}

\begin{remark}
For $\alpha>1$ we have $-2\alpha/D< -2/D$ and $(1-\alpha d)/(D-2)<-2/D$ so we get 
$S(X,t)\le C_d\tau^{\alpha-2/D}$ for $t$ large. The best inequality for $t$ large will be obtained for $d=2$. But in this case the constant can be improved by using directly lemma \ref{L:cased2}, as in the next proposition.
\end{remark}

\begin{proposition}
Let $\alpha>0$ and $t>2\pi$, then 
\begin{equation}\label{E:191218-1}
S(\tau^\alpha,t)\le 2A(\tau^{1/2}+2\tau^{\alpha-1/2}).
\end{equation}
\end{proposition}

\begin{proof}
Put $X=(t/2\pi)^\alpha=\tau^\alpha$, and $f(t)=-\frac{t}{2\pi}{\log x}$ we have 
$f'(t)=-\frac{t}{2\pi x}$ and $f''(x)=\frac{t}{2\pi x^2}$. The function $f''(x)$ satisfies the bounds
\[\lambda=\frac{t}{8\pi X^2}\le f''(x)\le\Lambda=\frac{t}{2\pi X^2},\qquad X<x\le 2X.\]
Applying Lemma 21 with $1<Y\le X$ yields
\begin{align*}
\Bigl|\sum_{X<n\le X+Y}e(f(n))\Bigr|&\le AY(\Lambda^2\lambda^{-1})^{1/2}+2A\lambda^{-1/2}
\le AX(\Lambda^2\lambda^{-1})^{1/2}+2A\lambda^{-1/2}\\
&=AX\Bigl(\frac{4t}{2\pi X^2}\Bigr)^{1/2}+2A\Bigl(\frac{t}{8\pi X^2}\Bigr)^{-1/2}=2A\tau^{1/2}+4A\tau^{\alpha-1/2}
\end{align*}
This inequality is also true for $0<Y\le 1$, since $2A\ge 5$, and $\tau\ge1$. Taking the supremum for $0<Y\le X$ we obtain
\eqref{E:191218-1}.
\end{proof}

\begin{remark}
It is easy to show that
\begin{equation} 
\begin{aligned}
\min_{d\ge2}\max(-2\alpha/D,&(1-\alpha d)/(D-2), -2/D)\\
&=\begin{cases}
-2/D & \text{for $\alpha\ge1$ with $d=2$},\\  
\displaystyle{\frac{1-\alpha d}{D-2}}& \text{for $\frac{2^{d-1}}{1+(d-1)2^{d-1}}<\alpha\le 
\frac{2^{d-2}}{1+(d-2)2^{d-2}}$ with $d$}.
\end{cases}
\end{aligned}
\end{equation}
Therefore, the exponent $\alpha$ determines which value of $d$ gives the least exponent. 
The next display shows the best values of $d$ for the given intervals of $\alpha$
\[\cdot \cdot, \frac{256}{2049}, d=9, \frac{128}{897}, d=8,  \frac{64}{385}, d=7;  \frac{32}{161}, d=6, \frac{16}{65}, d=5, \frac{8}{25}, d=4, \frac{4}{9}, d=3, \frac23, d=2. \]
\end{remark}

In applications, we divide the summation range into blocks; in some cases, the most important is the one with the largest exponent $\alpha$. In this situation, the following proposition is useful

\begin{proposition}\label{L:Sbound}
Let  $\alpha>0$ be a real number and $d\ge2$  an integer.  Assume that  $\alpha< \frac{2^{d-2}}{1+(d-2)2^{d-2}}$. For any $M\ge \widehat{B_d}$,  there exists a number $\tau_0=\tau_0(M)$ such that
\[S(\tau^\alpha,t)\le M \tau^{\alpha+\frac{1-\alpha d}{D-2}}\qquad \tau=t/2\pi\ge\tau_0.\]
The number $\tau_0$ is the least number such that  
\begin{equation}\label{E:t0}
\tau_0\ge d^{D/\alpha(D-2)},\quad \tau_0^{\frac{2}{D}+\frac{1-\alpha d}{D-2}}\ge\widehat{C_d}/M,\quad \tau_0^{\frac{2\alpha}{D}+\frac{1-\alpha d}{D-2}}\ge\widehat{A_d}/M.\end{equation}
\end{proposition}

\begin{proof}
First, we show that $-2/D\le -2\alpha/ D< (1-\alpha d)/(D-2)$.  By hypothesis 
$\alpha< \frac{2^{d-2}}{1+(d-2)2^{d-2}}\le 1$, and this proof the first inequality. For the second, notice that it is equivalent to $D-\alpha dD>-2\alpha D+4\alpha$, and therefore to
$D> \alpha((d-2)D+4)$ that by hypothesis is true. 

The exponents of $\tau_0$ in \eqref{E:t0} are positive. This implies that $\tau_0$ exists. 

Assume that $\tau>\tau_0$, with $\tau_0$ satisfying \eqref{E:t0}. Define $X=(t/2\pi)^\alpha$ since $\tau>\tau_0\ge1$, we have
\[X^{1-2/D}=\tau^{\alpha(D-2)/D}\ge\tau_0^{\alpha(D-2)/D}\ge d.\]
By Proposition \ref{P:30} we have 
\[S(X,t)\le\tau^\alpha\max\bigl\{\widehat{A_d}\tau^{-2\alpha/D}, \widehat{B_d}\tau^{\frac{1-\alpha d}{D-2}}, \widehat{C_d}\tau^{-2/D}\bigr\}\]
with $\widehat{A_d}$, $\widehat{B_d}$ and $\widehat{C_d}$ given by \eqref{constants}. 

Since the exponents in \eqref{E:t0} are positive, the inequalities there are also true for $\tau$ instead of $\tau_0$, and this implies that 
\[\max\bigl\{\widehat{A_d}\tau^{-2\alpha/D}, \widehat{B_d}\tau^{\frac{1-\alpha d}{D-2}}, \widehat{C_d}\tau^{-2/D}\bigr\}\le M\tau^{\frac{1-\alpha d}{D-2}},\] and the proof is complete.
\end{proof}

\begin{remark}
We can always take $M=\widehat{B_d}$, but we may obtain a too large value for $\tau_0$. 
The role of $M$ in Proposition \ref{L:Sbound} is  to allow a reasonable value for $\tau_0$, changing it to a larger coefficient. For example, the last two conditions in \eqref{E:t0} are automatically satisfied if we take $M=\max\{\widehat{A_d},\widehat{B_d},\widehat{C_d}\}$. 
\end{remark}

Also, the following proposition is useful to treat general zeta sums.
\begin{proposition}\label{P:partialsum}
Let $\sigma>0$,  $X>0$ and $t>1$ be  given, then 
\begin{equation}
S_\sigma(X,t):=\max_{X<Z\le 2X}\Bigl|\sum_{X<n\le Z}\frac{1}{n^{\sigma+it}}\Bigr|\le 
X^{-\sigma}S(X,t).
\end{equation}
\end{proposition}
\begin{proof}
For $y>0$ denote by $S(y)=\sum_{X<n\le y} n^{-it}$, so that $S(y)=0$ for $y\le X$.  By partial summation we have for any 
$X<Z\le 2X$
\[\sum_{X<n\le Z}\frac{1}{n^{\sigma+it}}=\frac{S(Z)}{Z^\sigma}+\sigma\int_X^Z S(y)y^{-\sigma-1}\,dy\]
Hence,
\[\Bigl|\sum_{X<n\le Z}\frac{1}{n^{\sigma+it}}\Bigr|\le \frac{S(X,t)}{Z^\sigma}+
\sigma S(X,t)\int_X^Z y^{-\sigma-1}\,dy=\frac{S(X,t)}{X^\sigma}.\]
Our result follows from taking supremum for $X<Z\le 2X$. 
\end{proof}

We finish by giving a classical example of an application to the Riemann zeta function.

\begin{proposition}
The following inequality holds 
\begin{equation}\label{E:zetabound}
\Bigl|\sum_{n\le (t/2\pi)^{1/2}}\frac{1}{n^{1/2+it}}\Bigr|\le
1.89725\; \tau^{1/6}\log\tau+9.89044\;\tau^{1/6},\qquad \text{for } \tau:=t/2\pi\ge648.
\end{equation}
\end{proposition}
\begin{proof}
We have 
\[\Bigl|\sum_{n\le \tau^\beta}\frac{1}{n^{1/2+it}}\Bigr|\le 
\sum_{n\le \tau^\beta}\frac{1}{n^{1/2}}\le \int_0^{\tau^\beta}u^{-1/2}\,du=2\tau^{\beta/2}.\]
Then we will have 
\begin{align*}
\Bigl|\sum_{n\le (t/2\pi)^{1/2}}\frac{1}{n^{1/2+it}}\Bigr|&\le
\Bigl|\sum_{n\le \tau^{1/3}}\frac{1}{n^{1/2+it}}\Bigr|+
\Bigl|\sum_{\tau^{1/3}<n\le\tau^{1/2}}\frac{1}{n^{1/2+it}}\Bigr|\\
&\le 2\, \tau^{1/6}+\sum_{k=1}^K \Bigl|\sum_{2^{-k}\tau^{1/2}<n\le 2^{-k+1}\tau^{1/2}}\frac{1}{n^{1/2+it}}\Bigr|,
\end{align*}
where $K$ is defined by the inequalities ${2^{-K}\tau^{1/2}<\tau^{1/3}\le 2^{-K+1}\tau^{1/2}}$, and in the sum for $n\le \tau^{1/3}$ we eliminate the terms appearing in the sum corresponding to $k=K$, the bound $2t^{1/6}$ is still valid for this sum.
Applying Proposition \ref{P:partialsum} yields 
\[\Bigl|\sum_{n\le (t/2\pi)^{1/2}}\frac{1}{n^{1/2+it}}\Bigr|\le 2\,\tau^{1/6}+
\sum_{k=1}^K(2^{-k}\tau^{1/2})^{-1/2} S(2^{-k}\tau^{1/2}, t).\]
Here $S(X_k,t)$ corresponds to the exponent $\alpha_k=\frac12-k\frac{\log2}{\log\tau}$
so that $\alpha_k<\frac12<\frac{2^{d-2}}{1+(d-2)2^{d-2}}$ for $d=3$. We want to apply Proposition \ref{L:Sbound}, to each of these $S(X_k,t)$, to this end we need $\tau$ satisfies \eqref{E:t0}, with $M=\widehat{B_d}$, that is,
\[\frac{D-2}{D}\Bigl(\frac12-k\frac{\log2}{\log\tau}\Bigr)\log\tau\ge\log d,\]\[
\Bigl\{\frac{2}{D}+\frac{1}{D-2}-\frac{d}{D-2}\Bigl(\frac12-k\frac{\log2}{\log\tau}\Bigr)\Bigr\}\log\tau\ge\log(\widehat{C_d}/\widehat{B_d}), \]
\[\Bigl\{\frac{1}{D-2}+\Bigl(\frac{2}{D}-\frac{d}{D-2}\Bigr)\Bigl(\frac12-k\frac{\log2}{\log\tau}\Bigr)\Bigr\}\log\tau\ge\log(\widehat{A_d}/\widehat{B_d}),\quad \text{for $1\le k\le K$}.\]
It is sufficient to check the first inequality for $k=K$ and the others for $k=1$. Since $d=3$ and $D=8$, the inequalities are equivalent to
\[\tau\ge2^{2K}3^{8/3},\quad 2^{3}\tau\ge(\widehat{C_3}/\widehat{B_3})^6,\quad 2^{6}\tau\ge (\widehat{A_3}/\widehat{B_3})^{24}\]
Since $2^K\le2\tau^{1/6}$ the first inequality follows from $\tau\ge2^{3}\cdot3^{4}$ and our conditions are 
\[\tau\ge2^{3}\cdot3^{4}=648,\quad \tau\ge (\widehat{C_3}/\widehat{B_3})^6/2^3=5.207\dots,\quad \tau\ge 127.537\dots.\]
Assuming this, we continue our bound 
\begin{align*}
\Bigl|\sum_{n\le (t/2\pi)^{1/2}}\frac{1}{n^{1/2+it}}\Bigr|&\le 2\,\tau^{1/6}+
\widehat{B_3}\sum_{k=1}^K(2^{-k}\tau^{1/2})^{-1/2} \tau^{\alpha_k+\frac{1-d\alpha_k}{D-2}} \\
&\le 2\,\tau^{1/6}+\widehat{B_3}\tau^{\frac16}\sum_{k=1}^K(2^{-k}\tau^{1/2})^{-1/2}
(2^{-k}\tau^{1/2})^{1-\frac12}\\
&=2\,\tau^{1/6}+\widehat{B_3}K\tau^{\frac16}\\
&\le 2\,\tau^{1/6}+\widehat{B_3}\tau^{\frac16}\Bigl(1+\frac16\frac{\log\tau}{\log2}\Bigr)
\end{align*}
Equation \eqref{E:zetabound} is obtained when we compute numerically the constants.
\end{proof}

\begin{remark}
Taking into account $\tau=t/2\pi$ equation \eqref{E:zetabound} is equivalent to 
\[\Bigl|\sum_{n\le (t/2\pi)^{1/2}}\frac{1}{n^{1/2+it}}\Bigr|\le 1.39668\; t^{1/6}\log t+4.71400\, t^{1/6}.\]
There are better bounds in the literature, but the bound has been obtained in a very straightforward way. 
\end{remark}

\section{Titchmarsh explicit Theorem}
The $d$-th derivative tests given usually are versions of Titchmarsh's Theorem 5.13, we will show that it is a consequence of van der Corput Theorem proving the following explicit version of it.
\begin{corollary}\label{TitchmarshTh}
 Let $f(x)$ be real and have continuous derivatives up to the order $d$-th, where $d\ge4$. Let $0<\lambda\le f^{(d)}(x)\le \Lambda$. Let $X$ and  $Y$ be real numbers  with $Y>d$ and $D=2^d$. Then,
\begin{equation}\label{E:corputgeneral}
\Bigl|\frac{1}{Y}\sum_{X<n\le X+Y}e^{2\pi i f(n)}\Bigr|\le C_4\max\Bigl\{
\Bigl(\frac{\Lambda}{\lambda}\Bigr)^{4/D}\lambda^{1/(D-2)}, Y^{-4/D}\lambda^{-1/(D-2)}\Bigr\}
\end{equation}
where $C_4\le 10.016$.
\end{corollary}
\begin{proof}
The trivial bound shows that the left-hand side of \eqref{E:corputgeneral} is  $E\le 1$. Hence we only have to proof our Theorem assuming that 
\[
\Bigl(\frac{\Lambda}{\lambda}\Bigr)^{4/D}\lambda^{1/(D-2)}<1,\qquad  Y^{-4/D}\lambda^{-1/(D-2)}<1.
\]
We note that these inequalities are equivalent to 
\begin{equation}\label{IneqInit}
\Lambda<\lambda^{1-\frac{D}{4(D-2)}},\quad \lambda^{-\frac{D}{4(D-2)}}<Y.
\end{equation}
We are in a position to apply Theorem \ref{T:dthtest} so that
\[E:=\Bigl|\frac{1}{Y}\sum_{X<n\le X+Y}e(f(n))\Bigr|\le
\max\Bigl\{A_d\Bigl(\frac{\Lambda}{\lambda Y}\Bigr)^{2/D}, B_d\Bigl(\frac{\Lambda^2}{\lambda}\Bigr)^{1/(D-2)},C_d(\lambda Y^d)^{-2/D}\Bigr\}.\]
The maximum value of $A_d$, $B_d$ and $C_d$ for $d\ge 4$ is $C_4=2(2AB^3)^{1/4}\le 10.016$. Therefore,
\[E\le C_4\max\Bigl\{\Bigl(\frac{\Lambda}{\lambda Y}\Bigr)^{2/D}, \Bigl(\frac{\Lambda^2}{\lambda}\Bigr)^{1/(D-2)},(\lambda Y^d)^{-2/D}\Bigr\}.\]
We end the proof showing that 
\begin{equation}\label{E:object}
\max\Bigl\{\Bigl(\frac{\Lambda}{\lambda Y}\Bigr)^{2/D}, \Bigl(\frac{\Lambda^2}{\lambda}\Bigr)^{1/(D-2)},(\lambda Y^d)^{-2/D}\Bigr\}\le \max\Bigl\{
\Bigl(\frac{\Lambda}{\lambda}\Bigr)^{4/D}\lambda^{1/(D-2)}, Y^{-4/D}\lambda^{-1/(D-2)}\Bigr\}.
\end{equation}

(a) The term $(\Lambda/\lambda Y)^{2/D}$ is less than $Y^{-4/D}\lambda^{-1/(D-2)}$ when $Y^{2/D}\le \lambda^{2/D-1/(D-2)}\Lambda^{-2/D}$. In the other case,
\[\Bigl(\frac{\Lambda}{\lambda Y}\Bigr)^{2/D}\le 
\Bigl(\frac{\Lambda}{\lambda }\Bigr)^{2/D} \lambda^{-2/D+1/(D-2)}\Lambda^{2/D}=\Bigl(\frac{\Lambda}{\lambda}\Bigr)^{4/D}\lambda^{1/(D-2)}.\]

(b) For the second term on the left-hand side of \eqref{E:object} we have
\[\Bigl(\frac{\Lambda^2}{\lambda}\Bigr)^{1/(D-2)}\le  \Bigl(\frac{\Lambda}{\lambda}\Bigr)^{4/D}\lambda^{1/(D-2)}.\]
This is true because it is equivalent to 
\[\Lambda^{\frac{2}{D-2}-\frac{4}{D}}\le \lambda^{\frac{2}{D-2}-\frac{4}{D}}\]
Since $\Lambda\ge\lambda$ this is equivalent to $\frac{2}{D-2}\le\frac{4}{D}$, that is true for $d\ge2$. 

(c) Only the case remains in which the maximum on the left-hand side of \eqref{E:object} is attained at $(\lambda Y^d)^{-2/D}$ so we may assume that 
\begin{equation}
\Bigl(\frac{\Lambda}{\lambda Y}\Bigr)^{2/D}<(\lambda Y^d)^{-2/D}, \quad \Bigl(\frac{\Lambda^2}{\lambda}\Bigr)^{1/(D-2)}<(\lambda Y^d)^{-2/D}.
\end{equation}
The first of these inequalities is equivalent to $Y<\Lambda^{-\frac{1}{d-1}}$. Since $Y\ge d\ge4$ it follows that $0<\lambda<\Lambda<1$. We divide the proof into two subcases according to which is the maximum of the two terms on the right-hand side of \eqref{E:object}.

(c1) Here, we assume that 
\[\Bigl(\frac{\Lambda}{\lambda}\Bigr)^{4/D}\lambda^{1/(D-2)}\le Y^{-4/D}\lambda^{-1/(D-2)}.\]
So, the maximum on the right-hand side of \eqref{E:object} is $Y^{-4/D}\lambda^{-1/(D-2)}$.
And we have to prove $(\lambda Y^d)^{-2/D}\le Y^{-4/D}\lambda^{-1/(D-2)}$. By contradiction assume that 
\[ (\lambda Y^d)^{-2/D}> Y^{-4/D}\lambda^{-1/(D-2)}\quad\text{i.~e.}\quad Y< \lambda^{-\frac{D-4}{2(D-2)(d-2)}}.\]
Then by \eqref{IneqInit} we have
\[\lambda^{-\frac{D}{4(D-2)}}<Y< \lambda^{-\frac{D-4}{2(D-2)(d-2)}}.\]
Since $\lambda<1$ this implies that 
\[\frac{D}{4(D-2)}<\frac{D-4}{2(D-2)(d-2)}.\]
And this is equivalent to saying that $dD+8<4D$. Therefore $d<4-8\times2^{-d}$, that contradicts our hypothesis $d\ge4$. 

(c2) Now, we assume that  
\[Y^{-4/D}\lambda^{-1/(D-2)}\le \Bigl(\frac{\Lambda}{\lambda}\Bigr)^{4/D}\lambda^{1/(D-2)}\quad\text{i.~e.}\quad \lambda^{1-\frac{D}{2(D-2)}}\Lambda^{-1}\le Y.\]
Hence, the  maximum on the right-hand side of \eqref{E:object} is $\left(\frac{\Lambda}{\lambda}\right)^{4/D}\lambda^{1/(D-2)}$. 
We have to prove that $(\lambda Y^d)^{-2/D}\le \left(\frac{\Lambda}{\lambda}\right)^{4/D}\lambda^{1/(D-2)}$. Arguing by contradiction, assume that 
\[\Bigl(\frac{\Lambda}{\lambda}\Bigr)^{4/D}\lambda^{1/(D-2)}<(\lambda Y^d)^{-2/D}\quad\text{i.~e.}\quad Y<\lambda^{\frac{D-4}{2d(D-2)}}\Lambda^{-\frac{2}{d}}.\]
Then, we obtain
\[\lambda^{1-\frac{D}{2(D-2)}}\Lambda^{-1}\le Y<\lambda^{\frac{D-4}{2d(D-2)}}\Lambda^{-\frac{2}{d}}.\]
It follows that 
\[\lambda^{1-\frac{D}{2(D-2)}}\Lambda^{-1+\frac{2}{d}}<\lambda^{\frac{D-4}{2d(D-2)}}\]
By \eqref{IneqInit} we have  $\Lambda<\lambda^{1-\frac{D}{4(D-2)}}$ and therefore 
\[\lambda^{\frac{2}{d}-\frac{D}{4(D-2)}-\frac{D}{2d(D-2)}}<\lambda^{\frac{D-4}{2d(D-2)}}.\]
Since $\lambda<1$, this gives us
\[\frac{D-4}{2d(D-2)}<\frac{2}{d}-\frac{D}{4(D-2)}-\frac{D}{2d(D-2)}\]
which is equivalent to $dD<4D-8$, which contradicts $d\ge 4$. 
\end{proof}

\end{document}